\documentclass[12pt,oneside,reqno]{amsart}
\usepackage{amsfonts}
\usepackage{amsmath}
\usepackage{amssymb}
\usepackage{geometry}

\DeclareFontFamily{U}{euc}{}
\DeclareFontShape{U}{euc}{m}{n}{<-6>eurm5<6-8>eurm7<8->eurm10}{}
\DeclareSymbolFont{AMSc}{U}{euc}{m}{n} 
\DeclareMathSymbol{\upgamma}{\mathord}{AMSc}{"D} 
\theoremstyle{plain}
\newtheorem{theorem}{Theorem}[section]

\newtheorem{proposition}[theorem]{Proposition}

\numberwithin{equation}{section}
\input{tcilatex}

\begin{document}
\title{Number identities and integer partitions}
\author{Craig Culbert}
\address{Department of Computer Science and Mathematical Sciences\\
Penn State University Harrisburg\\
Middletown, PA 17057}
\email{cwc15@psu.edu}
\date{May 23, 2018}

\begin{abstract}
Using a specific form of the triple product identity, polygonal number
identities are stated. Further number identities are examined that can be
considered identities related to modular sets of numbers. The identities can
be used to give results on integer partitions with parts from numbers in
modular arithmetic progression. This includes recursive formulas for the
number of partitions using these modular parts. The triple product identity
can derive further recursive formulas. Additionally, there is a recursive
formula for the related sum of divisors function. The specific triple
product identity provides a framework to examine all the identities and can
be used to define related theta functions.
\end{abstract}

\maketitle

\section{Number identities}

Euler's classic pentagonal number identity, Theorem 353 of \cite{hardywright}
or Corollary 1.7 of \cite{andrews3}, has exponents that are the general
pentagonal numbers indexed by the integers. There are identities that
involve triangular numbers and square numbers that can be found in Section
19.9 of \cite{hardywright}. Considering other polygonal numbers, Theorem 355
of \cite{hardywright} is an identity where the exponents are the general
heptagonal numbers. Additionally, two general polygonal number identities
attributed to Berger can be found in \cite{dickson} or \cite{berger}. All
these identities can be derived from Jacobi's triple product identity:
Theorem 352 of \cite{hardywright} or Theorem 2.8 of \cite{andrews3}. This
form of the triple product identity will be referred to as the traditional
form of the identity.

In Section 57\ of \cite{sylvester}, Sylvester references that Jacobi took
the traditional form of the triple product identity and through substitution
made the form more convenient for the study of number identities. This will
be the strategy employed in this article, but rather than use this form\ of
the identity, an alternative form of the triple product identity is used.
Denote the integers by $\mathbb{Z}$, the positive integers by $\mathbb{I}$
and the natural numbers including $0$ by $\mathbb{N}$. The variables $q$ and 
$z$ are complex numbers and consider the following form of the triple
product identity.

\begin{theorem}
For $\left\vert q\right\vert <1$ and $z\neq 0$,%
\begin{equation}
\prod\limits_{m\in \mathbb{I}}\left( 1-q^{m}\right) \!\!\left(
1+q^{m}z^{-1}\right) \!\!\left( 1+q^{m-1}z\right) =\sum\limits_{j\in \mathbb{%
Z}}q^{\frac{j^{2}-j}{2}}z^{j}  \label{triple}
\end{equation}
\end{theorem}

The identity is well know and can be proven by modifying the proof of
Theorem 2.8 of \cite{andrews3}. The identity can be found in a slightly
different form on page 285 of \cite{hardywright} and as Theorem 11 from
Chapter 8 of \cite{andrews1}. Following Jacobi, in the triple product
identity replace $q$ by $q^{k}$ and $z$ by $\pm q^{\ell }$. More generally, $%
z$ can be replaced by $q^{\ell }z$. As the eventual goal is to relate the
identities to integer partitions with positive integer parts, consider $k$ a
positive integer and $\ell $ a natural number restricted to $0\leq \ell \leq
k$. Each of the following infinite products is same for the pairs $\left(
k,\ell \right) $ and $\left( k,k-\ell \right) $. The equivalent series have
the same terms except they are indexed distinctly. Consequently, it is
possible to consider the restriction on $\ell $ of $0\leq \ell \leq \frac{k}{%
2}$.

\begin{align}
\prod\limits_{m\in \mathbb{I}}\left( 1-q^{km}\right) \!\!\left( 1+q^{km-\ell
}\right) \!\!\left( 1+q^{k\left( m-1\right) +\ell }\right) &
=\sum\limits_{j\in \mathbb{Z}}q^{\frac{k}{2}j\left( j-1\right) +\ell j}
\label{special one} \\
&  \notag \\
\prod\limits_{m\in \mathbb{I}}\left( 1-q^{km}\right) \!\!\left( 1-q^{km-\ell
}\right) \!\!\left( 1-q^{k\left( m-1\right) +\ell }\right) &
=\sum\limits_{j\in \mathbb{Z}}\left( -1\right) ^{j}q^{\frac{k}{2}j\left(
j-1\right) +\ell j}  \label{special two}
\end{align}%
For $\mathbf{\gamma }=\pm 1$, the products can be written using $q$%
-Pochhammer symbols as:%
\begin{equation*}
\left( q^{k}\,;q^{k}\right) _{\infty }\!\left( -\mathbf{\gamma }q^{k-\ell
}\,;q^{k}\right) _{\infty }\!\left( -\mathbf{\gamma }q^{\ell
}\,;q^{k}\right) _{\infty }=\left( q^{k},-\mathbf{\gamma }q^{k-\ell },-%
\mathbf{\gamma }q^{\ell }\,;q^{k}\right) _{\infty }
\end{equation*}

The two identities are similar to the special identities of section 19.9 of 
\cite{hardywright}. The second identity can be related naturally to integer
partitions and is similar to Corollary 2.9 of \cite{andrews3} where the $k$
is replaced by $2k+1$. It is possible now to consider polygonal number
identities for the general polygonal numbers, those indexed by $\mathbb{Z}$.
The classic polygonal numbers, those related to polygonal figures, are
indexed by $\mathbb{I}$ and can be found on page 1 of \cite{dickson}.
Consider $k$ a positive integer and let $\ell $ equal $1$, then the
following are the general polygonal number identities of Berger where $k+2$
is the number of sides of the polygon.

\begin{proposition}
For $\left\vert q\right\vert <1$,%
\begin{align*}
\prod\limits_{m\in \mathbb{I}}\left( 1-q^{km}\right) \!\!\left(
1+q^{km-1}\right) \!\!\left( 1+q^{k\left( m-1\right) +1}\right) &
=\sum\limits_{j\in \mathbb{Z}}q^{\frac{k}{2}j\left( j-1\right) +j} \\
& \\
\prod\limits_{m\in \mathbb{I}}\left( 1-q^{km}\right) \!\!\left(
1-q^{km-1}\right) \!\!\left( 1-q^{k\left( m-1\right) +1}\right) &
=\sum\limits_{j\in \mathbb{Z}}\left( -1\right) ^{j}q^{\frac{k}{2}j\left(
j-1\right) +j}
\end{align*}
\end{proposition}

These identities are found on page 31 of \cite{dickson} including a special
general hexagonal identity, or in \cite{berger}. Adding or subtracting a
pair of these identities gives further identities indexed by either even or
odd integers. In order to consider what classes of numbers are present in
the more general identities, specific sets of natural numbers are defined.

Let $k,\ell \in \mathbb{N}$, $k>0$ with $0\leq \ell \leq k$. Define $%
m_{k,\ell }(i)=k\left( i-1\right) +\ell $ for $i\in \mathbb{I}$. These
numbers constitute a \textit{modular arithmetic progression} or \textit{%
sequence} and $m_{k,\ell }(i)$ is the $i$\textit{-th number} in the
sequence. These numbers can be considered gnomons as defined on page 1 of 
\cite{dickson}. The sum of these gnomons leads to specific classes of
numbers, natural numbers greater than or equal to zero.

Consider the sum of these numbers in the progression, for $j\in \mathbb{I}$.%
\begin{equation*}
\sum\limits_{i=1}^{j}m_{k,\ell }(i)=\sum\limits_{i=1}^{j}k\left( i-1\right)
+\ell =\frac{k}{2}j\left( j-1\right) +\ell j=M_{k,\ell }(j)
\end{equation*}

The following is also true:%
\begin{equation*}
\sum\limits_{i=1}^{j}m_{k,k-\ell }(i)=\frac{k}{2}j\left( j+1\right) +-\ell j
\end{equation*}

These resulting sums form a class of numbers that can be considered figurate
in a very general sense. The set of these numbers $M_{k,\ell }(\mathbb{I}%
)=\left\{ M_{k,\ell }(j)\mid j\in \mathbb{I}\right\} $ are the \textit{%
modular figurate numbers} for parameters $k$ and $\ell $. $M_{k,\ell }(j)$
is the $j$\textit{-th number}. When indexed by the integers $\mathbb{Z}$,
these are the \textit{general modular figurate numbers}. Then $M_{k,\ell
}(-j)=M_{k,k-\ell }(j)$, $M_{k,\ell }(\mathbb{Z})=M_{k,k-\ell }(\mathbb{Z})$%
, and for a positive integer $c$, $M_{ck,c\ell }(j)=c\cdot \!M_{k,\ell }(j)$%
. All of the polygonal numbers and the general polygonal numbers can be
considered modular figurate numbers for $\ell =1$. Following the argument of
Berger for the polygonal numbers in \cite{berger}, for $i\neq j$, $M_{k,\ell
}(i)\neq M_{k,\ell }(j)$ unless $\ell =0$ with $i=1-j$, $\ell =k$ with $%
i=-1-j$, or $\ell =\frac{k}{2}$ with $i=-j$. These number forms are also
referenced in \cite{kolitsch}.

The modular figurate numbers are not the figurate numbers as defined on page
7 of \cite{dickson}. Instead, basic trapezoidal figures are related to the
form of these numbers, but the numbers are not all trapezoidal numbers. The
modular figurate numbers with the above restriction on $\ell $ can be
considered generalized trapezoidal numbers, the $k$-trapezoidal numbers
defined in \cite{jitphong} without the restriction on $\ell $.

Having defined this class of numbers, when (\ref{special one})\ and (\ref%
{special two}) are restricted to positive integers $k$ and natural numbers $%
\ell $ with $0\leq \ell \leq k$, these are identities of the modular
figurate numbers. Returning to the general setting of the triple product
identity, (\ref{triple}) can also be obtained by first considering the
following using Gaussian coefficients $\left[ \;\;\right] _{q}$.%
\begin{equation}
\prod\limits_{m=1}^{n}\left( 1+q^{m}z^{-1}\right) \!\!\left(
1+q^{m-1}z\right) =\sum\limits_{j=-n}^{n}\left[ 
\begin{array}{c}
2n \\ 
n+j%
\end{array}%
\right] _{q}q^{\frac{j^{2}-j}{2}}z^{j}  \label{hermite}
\end{equation}

The result can be proven using properties of the Gaussian coefficients and
induction. It is a modification of a result of Hermite found on page 49 of 
\cite{andrews3}. It is also possible to substitute in this equation $q^{k}$
and $\pm q^{\ell }$ to find a form with values $k$ and $\ell $. Of course
many of these identities are known and the traditional triple product
identity could be used to derive all the earlier identities. The use of the
present form of the triple product identity, (\ref{triple}), is to give a
specific symmetric framework to the study of the identities.

\section{Integer partitions}

In order to use (\ref{special two}) to study integer partitions with
distinct parts, the boundary conditions are first examined for both
identities and then excluded. Assuming $k$ is a positive integer with $0\leq
\ell \leq k$, one boundary case is when $\ell =0$ or by symmetry $\ell =k$.
The other case further assumes $k$ is an even positive integer and $\ell =%
\frac{k}{2}$. As the forms of these boundary cases are related to the
triangular and square numbers, their distinct properties can be derived from
the identities found in Section 19.9 of \cite{hardywright} and Corollary
2.10 of \cite{andrews3}.

For the case $\ell =0$, (\ref{special two}) is identically zero, but it is
still possible to combine the two identities to get specialized identities.
It is also possible to substitute $q^{k}$ into a triangular number identity,
(2.2.13) of \cite{andrews3}, to get a further identity for this boundary
case.

For the second boundary case when $k$ is an even positive integer and $\ell =%
\frac{k}{2}$, a substitution of $q^{\frac{k}{2}}$ into a square number
identity, (2.2.12) of \cite{andrews3}, gives a specialized identity for this
case. Another similar identity is the following.%
\begin{equation*}
\prod\limits_{m\in \mathbb{I}}\frac{\left( 1-q^{km}\right) \!\!\left(
1+q^{km-\frac{k}{2}}\right) }{\left( 1+q^{km}\right) \!\!\left( 1-q^{km-%
\frac{k}{2}}\right) }=\sum\limits_{j\in \mathbb{Z}}q^{\frac{k}{2}j^{2}}
\end{equation*}

In this case, where $k$ is a even positive integer, there is some
flexibility in combining (\ref{special one}) and (\ref{special two}) in
order to get identities indexed by odd or even integers, and to get
identities indexed by natural numbers or positive integers. One additional
identity is needed.%
\begin{equation*}
\prod\limits_{m\in \mathbb{I}}\left( 1+q^{km}\right) \!\!\left( 1+q^{km-%
\frac{k}{2}}\right) \!\!\left( 1-q^{km-\frac{k}{2}}\right) =1
\end{equation*}

Excluding the boundary conditions, the earlier identities give results about
integer partitions with distinct parts. Let $k$ and $\ell $ be positive
integers with $k\geq 3$ and $0<\ell <k$, $\ell \neq \frac{k}{2}$. Consider
sets of numbers from the modular arithmetic progression.%
\begin{align*}
\mathrm{J}_{k,\ell ,s}& =\left\{ m_{k,\ell }(i)\mid i\in \mathbb{I},1\leq
i\leq s\right\} \cup \left\{ m_{k,k-\ell }(i)\mid i\in \mathbb{I},1\leq
i\leq s\right\} \\
& \\
\mathrm{J}_{k,\ell }& =\left\{ m_{k,\ell }(i)\mid i\in \mathbb{I}\right\}
\cup \left\{ m_{k,k-\ell }(i)\mid i\in \mathbb{I}\right\}
\end{align*}%
Then $\mathrm{J}_{k,\ell }$ is the set of all positive integers congruent to
either $\pm \ell $ modulo $k$. It is also true that $\mathrm{J}_{k,\ell ,s}=%
\mathrm{J}_{k,k-\ell ,s}$ and $\mathrm{J}_{k,\ell }=\mathrm{J}_{k,k-\ell }$.
Let $\mathrm{J}$ be a set of positive integers and define $p_{_{\mathrm{dt}%
}}\!(n\,;\mathrm{J})$ as the number of partitions of the integer $n$ using
distinct parts only from the set $\mathrm{J}$. Denote the generating
function for $p_{_{\mathrm{dt}}}\!(n\,;\mathrm{J})$ as $f_{1}(q)$ as found
in \cite{andrews3}. The length of a partition is the number of parts and it
is possible to divide partitions into classes with an even or odd length.\
Define $r_{_{\mathrm{dt}}}\!(n\,;\mathrm{J})$ as the difference between the
number of the above partitions with even length and odd length, which has
integer entries and $r_{_{\mathrm{dt}}}\!(0\,;\mathrm{J})=1$. For negative
integers the number of partitions is always zero and denote the generating
function for $r_{_{\mathrm{dt}}}\!(n\,;\mathrm{J})$ as $g_{1}(q)$. In
general, generating functions for partition sequences that are the
difference between even and odd length can be found substituting $z=-1$ for
the length parameter in the two variable generating functions found in
Chapter 2 of \cite{andrews3}.

Consider now the modified result of Hermite, (\ref{hermite}). One
substitution gives an identity equal to the generating function of $p_{_{%
\mathrm{dt}}}\!(n\,;\mathrm{J}_{k,\ell ,s})$ and the other substitution
gives an identity equal to the generating function of $r_{_{\mathrm{dt}%
}}\!(n\,;\mathrm{J}_{k,\ell ,s})$. Both involve the modular figurate numbers 
$M_{k,\ell }(j)=\frac{k}{2}j\left( j-1\right) +\ell j$.%
\begin{align*}
\sum\limits_{n\in \mathbb{N}}p_{_{\mathrm{dt}}}\!(n\,;\mathrm{J}_{k,\ell
,s})q^{n}& =\sum\limits_{j=-s}^{s}\left[ 
\begin{array}{c}
2s \\ 
s+j%
\end{array}%
\right] _{q^{k}}q^{M_{k,\ell }(j)} \\
& \\
\sum\limits_{n\in \mathbb{N}}r_{_{\mathrm{dt}}}\!(n\,;\mathrm{J}_{k,\ell
,s})q^{n}& =\sum\limits_{j=-s}^{s}\left( -1\right) ^{j}\left[ 
\begin{array}{c}
2s \\ 
s+j%
\end{array}%
\right] _{q^{k}}q^{M_{k,\ell }(j)}
\end{align*}

It is now possible to write the above as a single formula, if some other
notation is introduced. Identify $\pm 1$ with $\pm $ and let $p_{_{\mathrm{dt%
}}}\!(n\,;\mathrm{J})=p_{_{\mathrm{dt}}}^{+}\!(n\,;\mathrm{J})$ and $r_{_{%
\mathrm{dt}}}\!(n\,;\mathrm{J})=p_{_{\mathrm{dt}}}^{-}\!(n\,;\mathrm{J})$.
This notation should not be confused with notation relating partitions to
even or odd permutations. It is possible to add or subtract the generating
functions to get expressions equivalent to the generating functions for
partitions on these sets of even or odd length.

For a positive integer $c$, define $c\mathrm{J}=\left\{ ci\mid i\in \mathrm{J%
}\right\} $. Define a further set of parts, $\overline{\mathrm{J}}_{k,\ell
}=k\mathbb{I}\cup \mathrm{J}_{k,\ell }$. Then $\overline{\mathrm{J}}_{k,\ell
}$ is the set of all positive integers congruent to either $0$ or $\pm \ell $
modulo $k$, and $\overline{\mathrm{J}}_{k,\ell }=\overline{\mathrm{J}}%
_{k,k-\ell }$. From this set of distinct parts (\ref{special two}) gives the
generating function for partition function $r_{_{\mathrm{dt}}}\!(n\,;%
\overline{\mathrm{J}}_{k,\ell })$.%
\begin{equation*}
g_{1}(q)=\sum\limits_{n\in \mathbb{N}}r_{_{\mathrm{dt}}}\!(n\,;\overline{%
\mathrm{J}}_{k,\ell })q^{n}=\sum\limits_{j\in \mathbb{Z}}\left( -1\right)
^{j}q^{\frac{k}{2}j\left( j-1\right) +\ell j}
\end{equation*}

Restating the identity as a result on integer partitions is this proposition
for $k\geq 3$, $0<\ell <k$ and $\ell \neq \frac{k}{2}$.

\begin{proposition}
Let $n$ be an integer.\label{rdistinct}%
\begin{equation*}
r_{_{\mathrm{dt}}}\!(n\,;\overline{\mathrm{J}}_{k,\ell })=%
\begin{cases}
\left( -1\right) ^{j} & \text{if }n=M_{k,\ell }(j)\text{,} \\ 
\;\;\;0 & \text{otherwise.}%
\end{cases}%
\end{equation*}
\end{proposition}

This result is similar to Theorem 2.11 of \cite{andrews3}. It appears that
this proposition is due to Sylvester as it is referenced in \cite{kolitsch}
and included there among combinatorial proofs of the triple product
identity. For the positive integers $\mathbb{I}=\overline{\mathrm{J}}_{3,1}$%
, Franklin's combinatorial proof of the partition result derives the
pentagonal number identity, Theorem 1.6 of \cite{andrews3}.

Euler used the values of $r_{_{\mathrm{dt}}}\!(n\,;\mathbb{I})$ to derive a
recursive formula for the partition function $p(n)$, Corollary 1.8 of \cite%
{andrews3}. A similar idea is now possible for unrestricted partitions on a
specific set of positive integers parts. Define the number of partitions of
an integer $n$ using only parts from a set $\mathrm{J}$ as $p(n\,;\mathrm{J}%
) $ with generating function $f(q)$ as in \cite{andrews3} and similarly
define $r(n\,;\mathrm{J})$ with generating function $g\left( q\right) $. The
generating functions for $r_{_{\mathrm{dt}}}\!(n\,;\mathrm{J})$ and $p(n\,;%
\mathrm{J})$ as infinite products are reciprocals, consequently the
generating functions are product inverses, $g_{1}(q)f(q)=1$. The values for $%
r_{_{\mathrm{dt}}}\!(n\,;\overline{\mathrm{J}}_{k,\ell })$ lead to a
recursive formula for $p(n\,;\overline{\mathrm{J}}_{k,\ell })$ as the number
of partitions for negative integers is zero and $p(0\,;\overline{\mathrm{J}}%
_{k,\ell })=1$.

\begin{theorem}
Let $n$ be a positive integer.%
\begin{gather*}
p(n\,;\overline{\mathrm{J}}_{k,\ell })=\sum\limits_{j\in \mathbb{Z-}\left\{
0\right\} }\left( -1\right) ^{j-1}p(n-M_{k,\ell }(j)\,;\overline{\mathrm{J}}%
_{k,\ell }) \\
\\
=\sum\limits_{j\in \mathbb{I}}\left( -1\right) ^{j-1}p(n-M_{k,\ell }(j)\,;%
\overline{\mathrm{J}}_{k,\ell })+\sum\limits_{j\in \mathbb{I}}\left(
-1\right) ^{j-1}p(n-M_{k,k-\ell }(j)\,;\overline{\mathrm{J}}_{k,\ell })
\end{gather*}
\end{theorem}

For a positive integer $c$, $p(cn\,;\overline{\mathrm{J}}_{ck,c\ell
})=p(cn\,;c\overline{\mathrm{J}}_{k,\ell })=p(n\,;\overline{\mathrm{J}}%
_{k,\ell })$ and the theorem gives a distinct recursive formula when $k$ and 
$\ell $ are relatively prime. Considering other sets of positive integer
parts, it is also possible to study the set of parts $k\mathbb{I}$, $k\geq 1$%
, using the results of Euler. If $\mathrm{I}_{k,\ell }$ is the set of
positive integers congruent to $\ell $ modulo $k$, substitutions into
Corollary 2.2 of \cite{andrews3} give results for the partition functions: $%
p(n)$, $r(n)$, $p_{_{\mathrm{dt}}}\!(n)$, $r_{_{\mathrm{dt}}}\!(n)$ and the
parts $\mathrm{I}_{k,\ell }$.

\section{Recursive formulas}

The previous recursive formula is derived from the triple product identity
and it is possible to identify other recursive formulas. Consider the
product of two generating functions equaling a third, $a(q)b(q)=c(q)$, then
there is a relationship between their coefficients:%
\begin{equation*}
\sum\limits_{k=0}^{n}a(k)b(n-k)=c(n)
\end{equation*}

If both $a(q)$ and $c(q)$ are generating functions with known formulas from
either (\ref{special one}) or (\ref{special two}), then a recursive formula
for $b(n)$\ can be derived by isolating the term $b(n-0)$. There are two
cases for the coefficients $c(n)$, nonzero or zero, and only the nonzero $%
a(k)$ summands are used. If only $a(q)$ has a known formula, then there is
an identity relating the sequence terms $b(n)$ and $c(n)$. In \cite{ono1},
the triple product identity is used to derive specific recursive formulas
and a similar procedure will be used focusing on examples that exclude the
boundary cases. By excluding examples that use triangular and square
numbers, the remaining modular figurate numbers are all indexed uniquely.

Let $\mathbf{\gamma }_{1}=\pm 1$, $\mathbf{\gamma }_{2}=\pm 1$ and consider
a function that is the quotient of two formulas derived from (\ref{triple}).
Assume neither product is identically zero. The function $H(q)$ converges
for $\left\vert q\right\vert <1$ and is a generating function indexed by $%
\mathbb{N}$.%
\begin{equation*}
H(q)=\frac{\left( q^{k_{2}},-\mathbf{\gamma }_{2}q^{k_{2}-\ell _{2}},-%
\mathbf{\gamma }_{2}q^{\ell _{2}}\,;q^{k_{2}}\right) _{\infty }}{\left(
q^{k_{1}},\mathbf{\gamma }_{1}q^{k_{1}-\ell _{1}},\mathbf{\gamma }%
_{1}q^{\ell _{1}}\,;q^{k_{1}}\right) _{\infty }}
\end{equation*}%
For the pairs $\left( k_{1},\ell _{1}\right) $ and $\left( k_{2},\ell
_{2}\right) $, let the associated modular figurate numbers be $M_{1}(\mathbb{%
Z})$ and $M_{2}(\mathbb{Z})$, respectively. Denote the sequence that has $%
H(q)$ as its generating function, $s^{\mathbf{\gamma }_{1},\mathbf{\gamma }%
_{2}}(n)$. Consider $H(q)$ multiplied by the product indexed by one and the
related generating functions. Using the values derived from either (\ref%
{special one}) or (\ref{special two}) gives a recursive formula for the
sequence $s^{\mathbf{\gamma }_{1},\mathbf{\gamma }_{2}}(n)$ as $s^{\mathbf{%
\gamma }_{1},\mathbf{\gamma }_{2}}(0)=1$ and consider $s^{\mathbf{\gamma }%
_{1},\mathbf{\gamma }_{2}}(n)=0$ for negative integers. For the following
theorem, assume that neither pair $\left( k_{1},\ell _{1}\right) $, $\left(
k_{2},\ell _{2}\right) $ is from a boundary case.

\begin{theorem}
Let $n$ be a positive integer. Let $\mathbf{\gamma }_{1}=\pm 1$ and $\mathbf{%
\gamma }_{2}=\pm 1$.%
\begin{equation*}
s^{\mathbf{\gamma }_{1},\mathbf{\gamma }_{2}}(n)=%
\begin{cases}
\sum\limits_{j\in \mathbb{Z-}\left\{ 0\right\} }-\left( -\mathbf{\gamma }%
_{1}\right) ^{j}s^{\mathbf{\gamma }_{1},\mathbf{\gamma }_{2}}(n-M_{1}(j))+%
\left( \mathbf{\gamma }_{2}\right) ^{i}\; & \text{if }n=M_{2}(i)\text{,} \\ 
&  \\ 
\sum\limits_{j\in \mathbb{Z-}\left\{ 0\right\} }-\left( -\mathbf{\gamma }%
_{1}\right) ^{j}s^{\mathbf{\gamma }_{1},\mathbf{\gamma }_{2}}(n-M_{1}(j)) & 
\text{otherwise.}%
\end{cases}%
\end{equation*}
\end{theorem}

If the function $H(q)$ involves a boundary case, that is the triangular or
square numbers, the argument is subtly different as evidenced by Theorem 1
or Theorem 3 of \cite{ono1}. Two consequences of the theorem for partition
sequences, $k\geq 3$ and $0<\ell <k$, $\ell \neq \frac{k}{2}$, now follow.
First, let $\left( k_{2},\ell _{2}\right) =\left( k,\ell \right) $ and $%
\left( k_{1},\ell _{1}\right) =\left( 3k,k\right) $ with $\mathbf{\gamma }%
_{1}=1$. Then $M_{1}(\mathbb{Z})=M_{3k,k}(\mathbb{Z})=kM_{3,1}(\mathbb{Z})=k%
\mathbb{I}$, $\mathbf{\gamma }=\pm 1$ and for function $H(q)$ the following
are equal. 
\begin{equation*}
\frac{\left( q^{k},-\mathbf{\gamma }q^{k-\ell },-\mathbf{\gamma }q^{\ell
}\,;q^{k}\right) _{\infty }}{\left( q^{3k},q^{2k},q^{k}\,;q^{3k}\right)
_{\infty }}=\frac{\left( q^{k}\,;q^{k}\right) _{\infty }\left( -\mathbf{%
\gamma }q^{k-\ell },-\mathbf{\gamma }q^{\ell }\,;q^{k}\right) _{\infty }}{%
\left( q^{k}\,;q^{k}\right) _{\infty }}
\end{equation*}%
Then $H(q)$ equals generating functions for partition functions involve the
sets $\mathrm{J}_{k,\ell }$.

\begin{equation*}
H(q)=\prod\limits_{m\in \mathbb{I}}\left( 1+\mathbf{\gamma }q^{km-\ell
}\right) \!\!\left( 1+\mathbf{\gamma }q^{k\left( m-1\right) +\ell }\right) 
\end{equation*}%
The recursive formulas that result use the notation for the general
pentagonal numbers, $\omega (j)=M_{3,1}(j)$.

\begin{proposition}
Let $n$ be a positive integer and $\mathbf{\gamma }=\pm 1$.%
\begin{equation*}
p_{_{\mathrm{dt}}}^{\mathbf{\gamma }}\!(n\,;\mathrm{J}_{k,\ell })=%
\begin{cases}
\sum\limits_{j\in \mathbb{Z-}\left\{ 0\right\} }\left( -1\right) ^{j-1}p_{_{%
\mathrm{dt}}}^{\mathbf{\gamma }}\!(n-k\omega (j)\,;\mathrm{J}_{k,\ell
})+\left( \mathbf{\gamma }\right) ^{i}\; & \text{if }n=M_{k,\ell }(i)\text{,}
\\ 
&  \\ 
\sum\limits_{j\in \mathbb{Z-}\left\{ 0\right\} }\left( -1\right) ^{j-1}p_{_{%
\mathrm{dt}}}^{\mathbf{\gamma }}\!(n-k\omega (j)\,;\mathrm{J}_{k,\ell }) & 
\text{otherwise.}%
\end{cases}%
\end{equation*}
\end{proposition}

By a similar procedure, consider $\mathbf{\gamma }_{2}=-1$ and $M_{2}(%
\mathbb{Z})=M_{3k_{1},k_{1}}(\mathbb{Z})=k_{1}M_{3,1}(\mathbb{Z})=k_{1}%
\mathbb{I}$. This leads to a recursive formula for unrestricted partitions
and the set $\mathrm{J}_{k,\ell }$.

\begin{proposition}
Let $n$ be a positive integer and $\mathbf{\gamma }=\pm 1$.%
\begin{equation*}
p^{\mathbf{\gamma }}(n\,;\mathrm{J}_{k,\ell })=%
\begin{cases}
\sum\limits_{j\in \mathbb{Z-}\left\{ 0\right\} }-\left( -\mathbf{\gamma }%
\right) ^{j}p^{\mathbf{\gamma }}(n-M_{k,\ell }(j)\,;\mathrm{J}_{k,\ell
})+\left( -1\right) ^{i}\; & \text{if }n=k\omega (i)\text{,} \\ 
&  \\ 
\sum\limits_{j\in \mathbb{Z-}\left\{ 0\right\} }-\left( -\mathbf{\gamma }%
\right) ^{j}p^{\mathbf{\gamma }}(n-M_{k,\ell }(j)\,;\mathrm{J}_{k,\ell }) & 
\text{otherwise.}%
\end{cases}%
\end{equation*}
\end{proposition}

Again, by adding or subtracting a pair of generating functions gives a
generating function for partitions of even and odd length. If instead, for $%
\mathbf{\gamma }_{1}=1,$ the relationship between the generating functions
is $H(q)=g_{1}(q)f(q)$, then the following identities are possible.

\begin{proposition}
Let $n$ be a natural number and $\mathbf{\gamma }=\pm 1$.%
\begin{eqnarray*}
p_{_{\mathrm{dt}}}^{\mathbf{\gamma }}\!(n\,;\mathrm{J}_{k,\ell })
&=&\sum\limits_{j\in \mathbb{Z}}\left( \mathbf{\gamma }\right)
^{j}p(n-M_{k,\ell }(j)\,;k\mathbb{I}) \\
p(n\,;\mathrm{J}_{k,\ell }) &=&\sum\limits_{j\in \mathbb{Z}}\left( -1\right)
^{j}p(n-k\omega (j)\,;\overline{\mathrm{J}}_{k,\ell })
\end{eqnarray*}
\end{proposition}

Theorem 5 of \cite{ono1} contains the generating function for partitions
excluding the parts congruent to $0$ modulo $d+1$. This is also equal to the
partitions where no part is repeated more than $d$ times, Corollary 1.3 of 
\cite{andrews3}. It is the later form that will be generalized. For a set of
positive integers $\mathrm{J}$, the number of partitions of an integer $n$
where no part appears more than $d$ times is denoted by $p_{_{\hat{d}}}(n\,;%
\mathrm{J})$, $d\geq 1$. In this way, $p_{_{\hat{1}}}(n\,;\mathrm{J})=p_{_{%
\mathrm{dt}}}\!(n\,;\mathrm{J})$. The generating function for $p_{_{\hat{d}%
}}(n\,;\mathrm{J})$ is found in \cite{andrews3}. Consider $\mathbf{\gamma }%
_{1}=1$, $\mathbf{\gamma }_{2}=-1$ and $\left( k_{2},\ell _{2}\right)
=\left( \left( d+1\right) k,\left( d+1\right) \ell \right) $, then $H(q)$
equals the following:

\begin{equation*}
H(q)=\frac{\left( q^{k_{2}},q^{k_{2}-\ell _{2}},q^{\ell
_{2}}\,;q^{k_{2}}\right) _{\infty }}{\left( q^{k},q^{k-\ell },q^{\ell
}\,;q^{k}\right) _{\infty }}=\prod\limits_{m\in \overline{\mathrm{J}}%
_{k,\ell }}\frac{1-q^{\left( d+1\right) m}}{1-q^{m}}
\end{equation*}

$H(q)$ is the generating function for $p_{_{\hat{d}}}(n\,;\overline{\mathrm{J%
}}_{k,\ell })$ which gives the following recursive formula.

\begin{proposition}
Let $n$ be a positive integer.%
\begin{equation*}
p_{_{\hat{d}}}(n\,;\overline{\mathrm{J}}_{k,\ell })=%
\begin{cases}
\sum\limits_{j\in \mathbb{Z-}\left\{ 0\right\} }\left( -1\right) ^{j-1}p_{_{%
\hat{d}}}(n-M_{k,\ell }(j)\,;\overline{\mathrm{J}}_{k,\ell })+\left(
-1\right) ^{i}\; & \text{if }n=\left( d+1\right) M_{k,\ell }(i)\text{,} \\ 
&  \\ 
\sum\limits_{j\in \mathbb{Z-}\left\{ 0\right\} }\left( -1\right) ^{j-1}p_{_{%
\hat{d}}}(n-M_{k,\ell }(j)\,;\overline{\mathrm{J}}_{k,\ell }) & \text{%
otherwise.}%
\end{cases}%
\end{equation*}
\end{proposition}

The case of $\overline{\mathrm{J}}_{3,1}=\mathbb{I}$ is Theorem 5 of \cite%
{ono1} and there is also an identity between partition values.

\begin{proposition}
Let $n$ be a natural number.%
\begin{equation*}
p_{_{\mathbf{\hat{\mathnormal{d}}}}}(n\,;\overline{\mathrm{J}}_{k,\ell
})=\sum\limits_{j\in \mathbb{Z}}\left( -1\right) ^{j}p(n-\left( d+1\right)
M_{k,\ell }(j)\,;\overline{\mathrm{J}}_{k,\ell })
\end{equation*}
\end{proposition}

Chapter 12 of \cite{rademacher} contains three relationships between
partition functions and the sum of divisors function, $\sigma (n)$. These
are a recursion for $p(n)$ using $\sigma (n)$, a recursion for $\sigma (n)$
using $r_{_{\mathrm{dt}}}\!(n\,;\mathbb{I})$ and the general pentagonal
numbers, and an identity for $\sigma (n)$ using $p(n)$ and the general
pentagonal numbers. It is possible to generalize these three concepts for
partitions in this article. Theorem 14.8 of \cite{apostal} provides formulas
to derive recursions involving $p(n\,;\mathrm{J})$ or $r_{_{\mathrm{dt}%
}}\!(n\,;\mathrm{J})$, and the following divisors function. Let $\mathrm{J}$
be a set of positive integers and $n\in \mathbb{I}$.%
\begin{equation*}
f_{\mathrm{J}}(n)=\sum\limits_{\substack{ d|n  \\ d\in \mathrm{J}}}d
\end{equation*}

For $n\in \mathbb{Z}$, $f_{\mathrm{J}}(n)$ is zero for negative integers and
zero. Denote the generating function for $f_{\mathrm{J}}(n)$ as $F(q)$. The
exposition for Theorem 14.8 of \cite{apostal} gives the relationships
between the generating functions, $qf^{\prime }(q)=f(q)F(q)$ and $%
qg_{1}^{\prime }(q)=-F(q)g_{1}(q)$, through logarithmic differentiation.
These relationships can be used to give a recursive formulas using $f_{%
\mathrm{J}}(n)$ for $p(n\,;\mathrm{J})$ and $r_{_{\mathrm{dt}}}\!(n\,;%
\mathrm{J})$, respectively. In order to the use the values for $r_{_{\mathrm{%
dt}}}\!(n\,;\overline{\mathrm{J}}_{k,\ell })$, Proposition \ref{rdistinct},
define $f_{k,\ell }(n)=f_{\mathrm{J}}(n)$ for the sets $\overline{\mathrm{J}}%
_{k,\ell }$. Such a divisors function is found in \cite{kim} using a
different notation. Theorem 14.8 of \cite{apostal} gives the following
equivalence.%
\begin{equation*}
n\cdot r_{_{\mathrm{dt}}}\!(n\,;\overline{\mathrm{J}}_{k,\ell })=-f_{k,\ell
}(n)-\sum\limits_{j=1}^{n-1}r_{_{\mathrm{dt}}}\!(j\,;\overline{\mathrm{J}}%
_{k,\ell })f_{k,\ell }(n-j)
\end{equation*}

Isolating $f_{k,\ell }(n)$, using the values for $r_{_{\mathrm{dt}}}\!(n\,;%
\overline{\mathrm{J}}_{k,\ell })$ and indexing over the nonzero integers
gives a finite recursive formula for this sum of divisors function for $%
k\geq 3$ and $0<\ell <k$, $\ell \neq \frac{k}{2}$.

\begin{theorem}
Let $n$ be a positive integer.%
\begin{equation*}
f_{k,\ell }(n)=%
\begin{cases}
\sum\limits_{j\in \mathbb{Z-}\left\{ 0\right\} }\left( -1\right)
^{j-1}f_{k,\ell }(n-M_{k,\ell }(j))+\left( -1\right) ^{i-1}M_{k,\ell }(i)\;
& \text{if }n=M_{k,\ell }(i)\text{,} \\ 
&  \\ 
\sum\limits_{j\in \mathbb{Z-}\left\{ 0\right\} }\left( -1\right)
^{j-1}f_{k,\ell }(n-M_{k,\ell }(j)) & \text{otherwise.}%
\end{cases}%
\end{equation*}
\end{theorem}

Chapter 12 of \cite{rademacher} provides the process to find an identity
that relates $f_{k,\ell }(n)$ and $p(n\,;\overline{\mathrm{J}}_{k,\ell })$.
For this remaining identity, $F(q)=-qg_{1}^{\prime }(q)f(q)$ is the
relationship between the generating functions. This results in an identity
of S. Kim found in \cite{kim}. There the identity was proven without the use
of Proposition \ref{rdistinct}, and the notation differs from the notation
in this article. As with the above recursions, it is possible to state this
identity of Kim using the modular figurate numbers $M_{k,\ell }(j)$.

\section{Further concepts}

Returning to the triple product identity (\ref{triple}), it is possible to
consider a theta function and related functions from the identity. Chapter
10 of \cite{rademacher} provides a guide to the motivation to define a theta
function. Consider for complex $\nu ,\tau $ with $\func{Im}(\tau )>0$.%
\begin{equation*}
\Theta (\nu \,;\tau )=\sum\limits_{n\in \mathbb{Z}}\exp \!\!\left( 2\pi 
\mathrm{i}\frac{n\left( n-1\right) }{2}\tau +2\pi \mathrm{i}n\nu \right)
\end{equation*}

The function $\Theta $ has the properties that $\Theta (\nu +\tau \,;\tau
)=\exp (-2\pi \mathrm{i}\nu )\Theta (\nu \,;\tau )$ and $\Theta (\nu
+1\,;\tau )=\Theta (\nu \,;\tau )$. With the substitutions $q=\exp (2\pi 
\mathrm{i}\tau )$ and $z=\exp (2\pi \mathrm{i}\nu )$ this gives the form in
the triple product identity above.%
\begin{equation*}
\Theta (z\,|\,q)=\sum\limits_{n\in \mathbb{Z}}q^{\frac{n\left( n-1\right) }{2%
}}z^{n}
\end{equation*}

It is then true that $\Theta (qz\,|\,q)=z^{-1}\Theta (z\,|\,q)$ and it is
possible to define auxiliary theta functions similar to the Jacobi theta
functions.%
\begin{align*}
\theta _{a}(z\,|\,q)& =\Theta (z\,|\,q) \\
\theta _{b}(z\,|\,q)& =\Theta (-z\,|\,q) \\
\theta _{c}(z\,|\,q)& =q^{-\frac{1}{8}}z^{\frac{1}{2}}\Theta (q^{\frac{1}{2}%
}z\,|\,q) \\
\theta _{d}(z\,|\,q)& =q^{-\frac{1}{8}}z^{\frac{1}{2}}\Theta (-q^{\frac{1}{2}%
}z\,|\,q)
\end{align*}

Each of these four auxiliary functions, when thought of as functions of $\nu 
$ and $\tau $, is a solution to a specific second order partial differential
equation. If the four functions are considered one class of theta functions,
then the substitutions $q$ replaced by $q^{k}$ and $z$ replaced by $q^{\ell
}z$ give another class of theta functions. The pair $\left( k,l\right) $
equal to $\left( 2,1\right) $ gives a form of the Jacobi theta functions. In
general, the functions $\theta _{a}$ and $\theta _{b}$ have infinite product
forms derived from (\ref{triple}) and product forms can also be determined
for $\theta _{c}$ and $\theta _{d}$.

Each class of theta functions potentially could be useful in the study of
the number of representations as sums of general polygonal numbers or
general modular figurate numbers. In Chapter 3 of \cite{mckeanmoll}, the
explanation of the Jacobi theta functions includes their use in study of the
number of representations as sums of squares. This involves identities of
the null values of the theta functions. The functions most similar to those
null values would be functions such as $\theta _{a}(\tau \,;2\tau )$.

The above theta functions are more closely related to the triangular
numbers. The use of theta functions in the study of the number of
representations as sums of triangular numbers can be found in \cite{ono2}.
In order to use such functions of $\tau $ to study sums of general polygonal
numbers or general modular figurate numbers, identities of these functions
would first need to be derived. That is to consider identities of the
functions either for a specific class of theta functions or in general
subject to the parameters $k$ and $\ell $.

\end{document}